\documentclass[11pt,a4paper]{article}
\usepackage{mycommands}

\newcommand{\inv}{\mathsf{inv}}
\newcommand{\td}{\mathsf{td}}

\newcommand{\Sn}[1][n]{\SSS_{#1}}
\newcommand{\Snm}{\Sn[n,m]}

\newcommand{\bss}{{\boldsymbol\sigma}}
\newcommand{\ssn}[1][n]{\bss_{#1}}
\newcommand{\ssk}{\ssn[k]}
\newcommand{\ssnm}[1][m]{\ssn[n,#1]}

\DeclareRobustCommand{\boldplus}{\mathord{\begin{tikzpicture}\draw[line width=0.3ex, x=1.5ex, y=1.5ex] (0.5,-0.05) -- (0.5,1.05)(-0.05,0.5) -- (1.05,0.5);\end{tikzpicture}}}
\DeclareRobustCommand{\boldtimes}{\mathord{\begin{tikzpicture}\draw[line width=0.3ex, x=1.5ex, y=1.5ex] (0.1,0.1) -- (0.9,0.9)(0.1,0.9) -- (0.9,0.1);\end{tikzpicture}}}

\newcommand{\In}[1][n]{\III_{#1}}

\newcommand{\Cn}[1][n]{\CCC_{#1}}

\newcommand{\ccn}[1][n]{{\boldsymbol\chi}_{#1}}

\newcommand{\dinv}{\rho_\inv}

\newcommand{\Nknml}[1][\ell]{N^{k,#1}_{n,m}}
\newcommand{\Nknmx}[1][\ell]{N^{k\+#1}_{n,m}}
\newcommand{\Nknmup}{\Nknmx[\nearrow]}
\newcommand{\Nknmdn}{\Nknmx[\searrow]}

\title{\textbf{Permutations with few inversions are locally uniform}}

\author{$\phantom{{}^\dagger}$David Bevan${}^\dagger$}

\hypersetup{
pdftitle={Permutations with few inversions are locally uniform},
pdfauthor={David Bevan},
pdfstartview={XYZ null null 1.00}
}

\date{}

\begin{document}
\maketitle

{\begin{NoHyper}
\let\thefootnote\relax
\end{NoHyper}}

{\begin{NoHyper}
\let\thefootnote\relax\footnotetext
{${}^\dagger$Department of Mathematics and Statistics, University of Strathclyde, Glasgow, Scotland.}
\end{NoHyper}}

\begin{abstract}
\noindent
We prove that permutations with few inversions exhibit a local--global dichotomy in the following sense.
Suppose $\bss$ is a permutation chosen uniformly at random from the set of all permutations of $[n]$ with exactly $m=m(n)\ll n^2$ inversions.
If $i<j$ are chosen uniformly at random from $[n]$, then $\bss(i)<\bss(j)$ asymptotically almost surely.
However, if $i$ and $j$ are chosen so that $j-i\ll m/n$, and $m \ll n^2/\log^2 n$, then $\liminfty\prob{\bss(i)<\bss(j)}=\half$.
Moreover, if $k=k(n)\ll \sqrt{m/n}$, then the restriction of $\bss$ to a random $k$-point interval is asymptotically uniformly distributed over $\Sn[k]$.
Thus, knowledge of the local structure of $\bss$ reveals nothing about its global form.
We establish that $k=\sqrt{m/n}$ is the threshold for the local uniformity of length $k$ subpermutations and that $j-i=m/n$ is the threshold for $\bss(i)\bss(j)$ to be as equally likely to be an inversion as not, and determine the behaviour in the critical windows.

\smallskip
\noindent
\textbf{\color{red}As pointed out by a referee, there are flaws in the proofs that do not seem easily rectifiable (see comments on pages~\pageref{ref1} and~\pageref{ref2}). So the results stated above have not been established.}
\end{abstract}

\section{Introduction}

We consider a permutation $\sigma$ of length $n$ (an \emph{$n$-permutation}) to be a linear ordering $\sigma(1)\ldots\sigma(n)$ of 
$[n]=\{1,\ldots,n\}$,
and identify $\sigma$ with its \emph{plot}, the set of points $\big\{(i,\sigma(i)) \::\: 1\leqslant i\leqslant n\big\}$ in the Euclidean plane.
For a very brief introduction to this perspective on permutations, see~\cite{BevanPPBasics}; for more extended expositions, see either B\'ona~\cite{BonaBook} or Kitaev~\cite{KitaevBook}; for an extensive recent survey, see Vatter~\cite{VatterSurvey}.
We use $\Sn$ for the set of all $n$-permutations and
$|\sigma|$ to denote the length of~$\sigma$.

An \emph{inversion} in $\sigma$ is a pair $i,j$ of indices such that $i<j$ and $\sigma(i)>\sigma(j)$, or equivalently two points in the plot of $\sigma$, one to the northwest of the other.
The number of inversions in $\sigma$ is denoted $\inv(\sigma)$,
and we use
$\Snm = \{ \sigma\in \Sn : \inv(\sigma) = m \}$ for the set of all $n$-permutations with exactly $m$ inversions.
An $n$-permutation can have at most $\binom{n}{2}$ inversions; the \emph{inversion density} of $\sigma$ is defined to be the ratio
$\dinv(\sigma) = \inv(\sigma)/\binom{|\sigma|}{2}$.

We are interested in studying the properties of a typical large $n$-permutation with $m=m(n)$ inversions, for a given function $m(n)$, and we use
$\ssnm$ to denote a permutation chosen uniformly at random from $\Snm$.
Thus, $\ssnm$ can be seen as the natural analogue for permutations of the \Erdos--R\'enyi random graph $\mathbf{G}_{n,m}$~\cite{ER1960}.
The only prior work on $\ssnm$
of which we are aware
is that of Acan and Pittel~\cite{AP2013}, who establish a sharp threshold for connectivity (that is, sum indecomposability) at $m=(6/\pi^2)n\log n$.

\together3
In this paper, our focus is on the \emph{local} structure of $\ssnm$ when
$m$ grows superlinearly but subquadratically with $n$, that is\footnote{We write $f(n)\ll g(n)$ or $g(n)\gg f(n)$ if $\liminftyt f(n)/g(n) = 0$, and write $f(n)\sim g(n)$ if $\liminftyt f(n)/g(n) = 1$.} when $n\ll m\ll n^2$.
We call such (random) permutations \mbox{\emph{semi-sparse}}.

\begin{figure}[t] 
  \centering
  \setlength{\fboxsep}{0pt}
  \fbox{\includegraphics[width=1.8in]{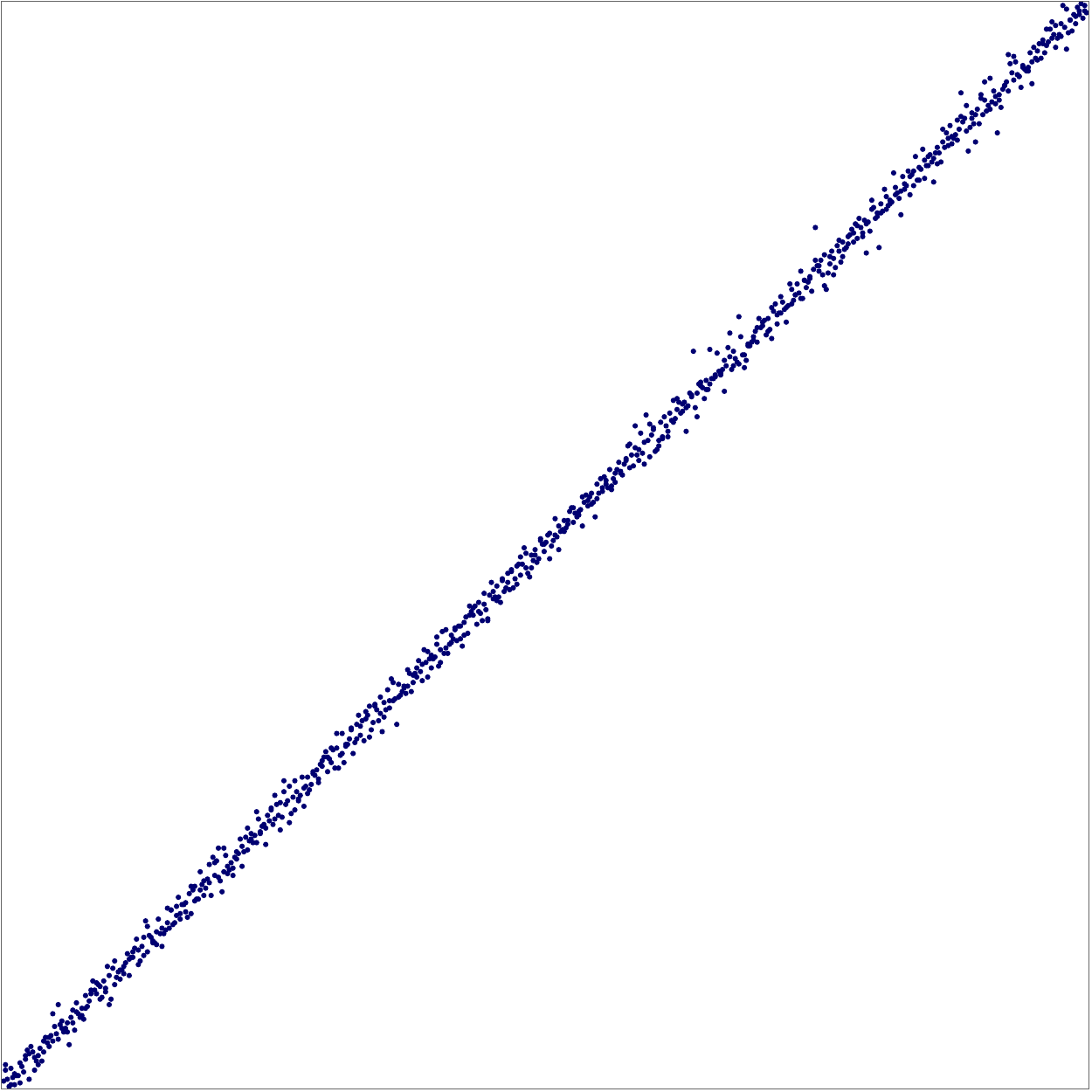}}
  \caption{The plot of a randomly selected permutation on 825 points with inversion density 0.01.}
  \label{figPerm}
\end{figure}

Almost all the points of a semi-sparse permutation are close to the main diagonal in the following sense.
The \emph{absolute displacement} of the $j$th point of $\sigma$ is $d_j(\sigma)=|\sigma(j)-j|$, its vertical distance from the main diagonal.
Knuth~\cite{Knuth1973} defined the \emph{total displacement} $\td(\sigma) = \sum_{j=1}^{|\sigma|} d_j(\sigma)$ as a natural measure of how close $\sigma$ is to the identity.
Subsequently, Diaconis and Graham~\cite{DG1977} proved that total displacement and number of inversions are related by the following double inequality:
$\inv(\sigma)\leqs\td(\sigma)\leqs2\+\inv(\sigma)$ for any permutation $\sigma$.
Consequently, if $\ssnm$ is semi-sparse,
\[
\liminfty \prob{d_{j_n}(\ssnm) \gg m/n} \;=\; 0 ,
\]
for any sequence of positive integers $(j_n)$ with $j_n\leqs n$.
Indeed, the analytic limit $\liminftyt\ssnm$ of semi-sparse permutations is the \emph{permuton}\footnote{A permuton is a probability measure on the unit square with uniform marginals.} whose support is the main diagonal
(see~\cite{KKRW2020} and references therein for information about permutons).
See Figure~\ref{figPerm} for an illustration of a permutation with small inversion density chosen uniformly from $\Sn[825,3399]$.

Unsurprisingly, then, if we pick two points randomly from a semi-sparse permutation $\ssnm$, then {asymptotically almost surely}\footnote{A property $Q=Q(n)$ holds {asymptotically almost surely} if $\liminftyt\prob{Q} = 1$.} they do not form an inversion:
\[
\text{If ~} i<j\text{, ~then ~} \prob{\ssnm(i)>\ssnm(j)} \;=\; \dinv(\ssnm) \;=\; m/\tbinom{n}{2} \;\to\; 0 \text{~ as ~} n\to\infty .
\]
However, in stark contrast, the \emph{local} structure is very different.
Our primary result is that if we ``zoom in'' far enough, then {locally} a semi-sparse permutation is \emph{uniform}, the restriction of $\ssnm$ to a sufficiently small interval being asymptotically uniformly distributed. The local structure of $\ssnm$ thus reveals nothing about its global form.
As an illustration of this phenomenon, even for relatively small $n$, in the permutation in Figure~\ref{figPerm} more than 47\% of the pairs of \emph{adjacent} points form inversions (that is,  descents).
For a rather different perspective on the independence of local and global structure in permutations, see~\cite{BP2020} and~\cite{Bevan2022}.

\together4
\subsection{Results}


To state our results, we use the following notation and phraseology:
$\sigma[i,j]$ 
denotes the sequence $\sigma(i)\sigma(i+1)\ldots\sigma(j)$,
the restriction of $\sigma$ to the interval $[i,j]$.
We say that $\sigma[i,j]$ \emph{forms}~$\tau$ if the terms of $\sigma[i,j]$ are in the same relative order as those of the permutation $\tau$,
and we say that
$\tau$ \emph{occurs at position $j$ in~$\sigma$} if $\sigma[j,j+|\tau|-1]$ forms $\tau$.
In this context, the (\emph{consecutive}) \emph{subpermutation} $\tau$ is called a (\emph{consecutive}) \emph{pattern}.

Our first theorem establishes uniformity at a sufficiently small scale.

\begin{thmO}\label{thmPatUnif}
  Suppose $n\ll m \ll n^2/\log^2n$ and $k \ll \sqrt{m/n}$.
  Then, for any sequence of permutations $(\tau_n)$ with
  $|\tau_n|=k$
  and
  any sequence of positive integers $(j_n)$ with $j_n\leqs n+1-k$,
  \[
  \prob{\text{$\tau_n$ occurs at position $j_n$ in $\ssnm$}}
  \;\sim\;
  \frac{1}{k!}
  .
  \]
\end{thmO}

If we ``zoom out'' a little, then we have the following behaviour in the critical window, where the probability that $\ssnm$ locally looks like a permutation $\tau$ depends on $\tau$'s {inversion density}.

\begin{thmO}\label{thmPatCrit}
  Suppose $n\ll m \ll n^2/\log^2n$ and $k\sim \alpha\sqrt{m/n}$ for some $\alpha>0$.
  Fix $\rho\in[0,1]$.
  Then, for any sequence of permutations $(\tau_n)$ with $|\tau_n|=k$ and
  $\dinv(\tau_n)\sim\rho$,
  and
  any sequence of positive integers $(j_n)$ with $j_n\leqs n+1-k$,
  \[
  \prob{\text{$\tau_n$ occurs at position $j_n$ in $\ssnm$}}
  \;\sim\;
  { e^{(1-2\rho)\alpha^2/4} }
  \frac {1}
  { k! }
  .
  \]
\end{thmO}

If our window is a bit wider,
then any subpermutation with sufficient inversion density almost never occurs.

\begin{thmO}\label{thmPat0}
  Suppose $n\ll m \ll n^2/\log^2n$ and $k\gg \sqrt{m/n}$.
  Suppose $m/nk^2\ll\rho\leqs1$.
  Then, for any sequence of permutations $(\tau_n)$ with $|\tau_n|=k$ and
  $\dinv(\tau_n)\sim\rho$,
  and
  any sequence of positive integers $(j_n)$ with $j_n\leqs n+1-k$,
  \[
  \prob{\text{$\tau_n$ occurs at position $j_n$ in $\ssnm$}}
  \;\ll\;
  \prob{\text{$\ssnm[m][j_n,j_n+k-1]$ is increasing}}
  .
  \]
\end{thmO}

These first three theorems thus reveal the existence of a
threshold at $k=\sqrt{m/n}$ for 
consecutive $k$-sub\-permutations of semi-sparse $\ssnm$ to be uniformly distributed.

If we turn our attention to pairs of points, the threshold for $i,j$ being an inversion in $\ssnm$ is at the larger scale of $j-i=m/n$.
Below this, $i,j$ is as likely to be an inversion as not, whereas above it, $i,j$ is almost never an inversion.

\begin{thmO}\label{thmInvUnif}
  Suppose $n\ll m \ll n^2/\log^2n$ and $(j_n)$ is any sequence of positive integers such that $j_n\leqs n-k$.
  Then,
  \[
  \liminfty\prob{\ssnm(j_n) > \ssnm(j_n+k)}
  \;=\;
  \begin{cases}
  \thalf & \text{if $k\ll {m/n}$,} \\[2pt]
  0 & \text{if $k\gg {m/n}$.}
  \end{cases}
  \]
\end{thmO}

To conclude, we determine for $j-i$ in the critical window the exact asymptotic probability of $i,j$ forming an inversion.

\begin{thmO}\label{thmInvCrit}
  Suppose $n\ll m \ll n^2/\log^2n$ and $k\sim \alpha {m/n}$ for some $\alpha>0$.
  Then, for any sequence of positive integers $(j_n)$ with $j_n\leqs n-k$,
  \[
  \liminfty\prob{\ssnm(j_n) > \ssnm(j_n+k)}
  \;=\;
  \frac{e^\alpha(\alpha-1)+1}{(e^\alpha-1)^2}
  \;<\;
  \half
  .
  \]
\end{thmO}

\subsection{Methodology}

The four essential components of our approach are as follows.
Firstly, we establish a position independence result (Proposition~\ref{propShiftInvariance}):
For any $n$ and $m$, the random permutation $\ssnm$ ``looks the same'' in any two intervals of the same length.
This means that we need only consider the structure of the first $k$ points of~$\ssnm$.
Secondly, we represent permutations by their inversion sequences.
Thirdly, we
establish conditions under which we can
use weak compositions to
approximate large suffixes of inversion sequences of random semi-sparse permutations (Corollary~\ref{corWeakComp}).
And, finally, we make use of a tripartition of
these inversion sequences to enable us to apply these approximations to asymptotically enumerate certain classes of permutations (Proposition~\ref{propNknml} and Proposition~\ref{propNComp}) which we use to yield our results.

In Section~\ref{sectFoundations}, we develop this framework, proving the position independence of consecutive patterns and establishing how to count semi-sparse permutations by approximating inversion sequences with weak compositions.
Section~\ref{sectUnif} then contains the proofs of Theorems~\ref{thmPatUnif}, \ref{thmPatCrit} and~\ref{thmPat0}, establishing the threshold for local uniformity.
Finally, in Section~\ref{sectInv}, we prove Theorems~\ref{thmInvUnif} and~\ref{thmInvCrit}, which establish the threshold for inversions.

The approximations we use are only valid when $m \ll n^2/\log^2n$.
(See the comment after the proof of Proposition~\ref{propNknml}.)
For faster-growing $m$, a different approach is needed.
Thus, the following questions remain open.

\begin{questO}\label{questUnif}
  Suppose $m=\Omega\big(n^2/\log^2n\big)$.
  How slowly does $k$ need to grow so that
  for any sequence of permutations $(\tau_n)$ with
  $|\tau_n|=k$
  and
  any sequence of positive integers $(j_n)$ with $j_n\leqs n+1-k$,
  \[
  \prob{\text{$\tau_n$ occurs at position $j_n$ in $\ssnm$}}
  \;\sim\;
  \frac{1}{k!}
  ?
  \]
\end{questO}

\begin{questO}\label{questInv}
  Suppose $m=\Omega\big(n^2/\log^2n\big)$.
  How slowly does $k$ need to grow so that
  for any sequence of positive integers $(j_n)$ with $j_n\leqs n-k$,
  \[
  \liminfty\prob{\ssnm(j_n) > \ssnm(j_n+k)}
  \;=\;
  \thalf
  ?
  \]
\end{questO}
It seems likely that techniques suitable for answering these questions for the remainder of the semi-sparse range ($m\ll n^2$)
may well not be applicable to dense permutations, when
$m\sim\rho\binom{n}{2}$.
The only simple solution appears to be for Question~\ref{questInv} in the dense case when $\rho=\half$, when the probability of $j,j+k$ forming an inversion is $\nhalf$ for all $k<n$.




\section{Foundations}\label{sectFoundations}

In this section, we establish the basic framework we use to prove our results: the position independence of subpermutations and the
asymptotic enumeration of inversion sequences by approximation using weak compositions.


\subsection{Position independence} 

Our first observation is that the distribution of any consecutive pattern in $\ssnm$ is independent of its position. This holds for any given $n$ and $m$.
As a consequence, in subsequent arguments, we need only consider the occurrence of patterns at position 1 in $\ssnm$.

\begin{propO}
\label{propShiftInvariance}
  For any permutation $\tau\in \Sn[k]$
  and any positive $i,j\leqs n+1-k$,
  \[
  \prob{\text{$\tau$ occurs at position $i$ in $\ssnm$}} \:=\: \prob{\text{$\tau$ occurs at position $j$ in $\ssnm$}}
  .
  \]
\end{propO}
This result follows from the existence of an operation that removes the last point from a permutation and adds a new first point in such a way as to preserve the number of inversions. This operation shifts 
patterns rightwards.
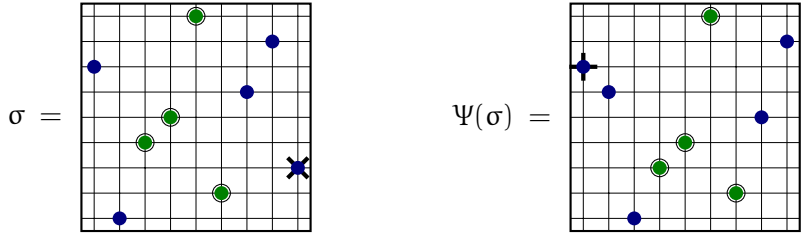
\begin{figure}[t] 
  \centering
\begin{tikzpicture}[scale=0.335]
    \draw[ultra thick] (8.6,2.6) -- (9.4,3.4);
    \draw[ultra thick] (9.4,2.6) -- (8.6,3.4);
    \plotpermgrid  [blue!50!black]{9}{7,1,0,0,0,0,6,8,3}
    \plotpermnobox[green!50!black]{9}{0,0,4,5,9,2}
    \circpt{3}{4}
    \circpt{4}{5}
    \circpt{5}{9}
    \circpt{6}{2}
    \node[left] at (.7,5) {$\sigma \:=\: {}$};
\end{tikzpicture}
$\qquad\qquad$
\begin{tikzpicture}[scale=0.335]
    \draw[ultra thick] (0.45,7) -- (1.55,7);
    \draw[ultra thick] (1,6.45) -- (1,7.55);
    \plotpermgrid  [blue!50!black]{9}{7,6,1,0,0,0,0,5,8}
    \plotpermnobox[green!50!black]{9}{0,0,0,3,4,9,2}
    \circpt{4}{3}
    \circpt{5}{4}
    \circpt{6}{9}
    \circpt{7}{2}
    \node[left] at (.7,5) {$\Psi(\sigma) \:=\: {}$};
\end{tikzpicture}
  \caption{The bijection used in the proof of Proposition~\ref{propShiftInvariance}:
  the point marked $\boldtimes$ is replaced by that marked $\boldplus$;
  the pattern 2341 occurs at position~3 in $\sigma$ and at position~4 in $\Psi(\sigma)$}
  \label{figShiftInvariance}
\end{figure}
\begin{proof}
  As illustrated in Figure~\ref{figShiftInvariance},
  let $\Psi:\Snm\to\Snm$ be defined by
  \[
  \Psi(\sigma) \;=\;
  \Psi(\sigma_1\sigma_2\ldots\sigma_n) \;=\; \sigma' \;=\; \sigma'_0\sigma'_1\ldots\sigma'_{n-1} ,
  \]
  where 
  $\sigma'_0=n+1-\sigma_n$, and for $1\leqs i<n$,
  $$
  \sigma'_i \;=\;
  \begin{cases}
  \sigma_i+1 , & \text{~if~ $\sigma'_0\leqs\sigma_i<\sigma_n$,} \\
  \sigma_i-1 , & \text{~if~ $\sigma_n<\sigma_i\leqs \sigma'_0$,} \\
  \sigma_i   , & \text{~otherwise.}
  \end{cases}
  $$
Note that $\sigma_n$ contributes $n-\sigma_n$ inversions to $\sigma$, and
$\sigma'_0$ contributes the same number of inversions to $\sigma'$.
For $0< i< n$, the point $\sigma'_i$ contributes the same number of inversions to $\sigma'$ as
$\sigma_i$ does to $\sigma$. So $\inv(\sigma')=\inv(\sigma)$.
Since $\Psi$ preserves length and has a well-defined inverse, it is a bijection on $\Snm$.

If $\tau\in \Sn[k]$ occurs at position $j\leqs n-k$ in $\sigma$, then $\tau$ occurs at position $j+1$ in $\Psi(\sigma)$.
Hence, if $1\leqs i,j \leqs n+1-k$, then $\tau$ occurs at position $i$ in $\sigma$ if and only if $\tau$ occurs at position $j$ in~$\Psi^{j-i}(\sigma)$.
\end{proof}

\subsection{Inversion sequences and weak compositions}

\begin{figure}[t] 
  \centering
  \begin{tikzpicture}[scale=0.31]
    \plotpermgrid[blue!50!black]{9}{7,3,5,8,4,6,1,9,2}
    \foreach \y [count=\x] in {0,1,1,0,3,2,6,0,7}{
      \node at (\x,-.4) {\small\textsf{\y}};
      \node at (\x+16,-.4) {\small\textsf{\y}};
      \draw[] (\x+16-.375,-10.45+\x+10) -- (\x+16-.375,-9.5+10) -- (\x+16+.375,-9.5+10) -- (\x+16+.375,-10.45+\x+10);
      \ifnum0=\y {} \else {
        \foreach \z in {1,...,\y} \fill[radius=0.3,blue!50!red!50!black] (\x+16,\z) circle;
      } \fi
    }
  \end{tikzpicture}
  \caption{A permutation with 20 inversions and its inversion sequence}
  \label{figInvSeq}
\end{figure}
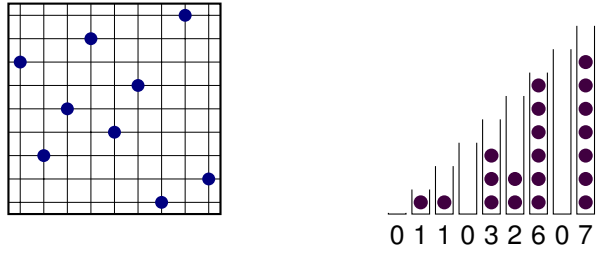

Key to our analysis is the representation of permutations as \emph{inversion sequences}.
The {inversion sequence} of an $n$-permutation $\sigma$ is $(e_j)_{j=1}^n$, where $e_j=\big|\{i: i<j \text{~and~} \sigma_i>\sigma_j \}\big|$
is the number of inversions of $\sigma$ whose right end is at position $j$.
Note that $\inv(\sigma)=\sum_j e_j$.

Clearly, for each $j$, we have $e_j<j$, and in fact sequences satisfying this condition whose sum equals~$m$ are in bijection with $n$-permutations having $m$ inversions.
Each $e_j$ can 
be considered to be the number of \emph{balls} in a \emph{box} whose capacity is $j-1$. This is illustrated in Figure~\ref{figInvSeq}.
Rather than working directly with permutations, we investigate the properties of this balls-in-boxes model when the number of balls ($m$) is superlinear but subquadratic in the number of boxes~($n$).

This analysis is aided by using (unrestricted) \emph{weak compositions} to approximate inversion sequences.
A \emph{weak $t$-composition of $s$} is a sequence of $t$ non-negative integers whose sum is~$s$.
In terms of our balls-in-boxes model, we have $s$ balls in $t$ boxes, each of whose capacity is \emph{unlimited}.
We use $\Cn[t,s]$ to denote the set of all weak $t$-compositions of $s$.
Clearly, the number of such compositions is given by $|\Cn[t,s]| = \binom{s+t-1}{s}$.

The foundation for our approximation is the fact that asymptotically almost surely, no term in a semi-sparse weak $t$-composition of $s$ exceeds
its expected value of $s/t$ by a factor
greater than $\log t$.
A specific instance of this proposition is proved by Acan and Pittel in~\cite{AP2013}; we generalise their proof.

\begin{propO}\label{propWeakComp}
  Suppose $t\ll s\ll t^2/\log t$ and $\ccn[t]$ is chosen uniformly at random from $\Cn[t,s]$.
  Then, for any $\veps>0$ and sufficiently large $t$, 
  \[
 \prob{\text{some term of $\ccn[t]$ is at least $(1+\veps)\tfrac{s}{t}\log t$}} \;\leqs\; t^{-\veps/2} .
  \]
\end{propO}
\begin{proof}
  Let $L_r$ be the number of 
  terms of $\ccn[t]$ whose value is at least $r$. By the first moment method and linearity of expectation, 
  $\prob{L_r > 0} \,\leqs\, \expec{L_r} \,=\, t\+\prob{x_1 \geqs r}$, where $x_1$ is the first term of $\ccn[t]$.

  Now, for any positive $r$,
  \[
   \prob{x_1 \geqs r}
   \;=\; \frac{\sum_{j=0}^{s-r}|\Cn[t-1,j]|}{|\Cn[t,s]|}
   \;=\; \frac{s!\+(s+t-r-1)!}{(s-r)!\+(s+t-1)!}
   \;\leqs\; \left(\frac{s}{s+t-r}\right)^{\!r}
   \;=\; \left(\!1-\frac{t-r}{s+t-r}\right)^{\!r}
  .
  \]
  Thus, $\prob{L_r > 0} \,\leqs\, t\!\left(1-\frac{t-r}{s+t-r}\right)^{\!r} \,\leqs\, t\exp\!\left(-\frac{r(t-r)}{s+t-r}\right)$.

  Suppose $r=\alpha\frac{s}{t}\log t$, 
  where $\alpha=1+\veps$.
  By the bound on $s$, we have $r<t$ for sufficiently large~$t$.
  Rearrangement then yields
  \[
  \frac{r(t-r)}{s+t-r}
  \;=\; \left(\frac{1-\frac{\alpha s \log t}{t^2}}{1+\frac{t}{s}-\frac{\alpha\log t}{t}}\right) \alpha\log t
  \;\geqslant\; \big(\alpha-\tfrac{\veps}{2}\big)\log t
  \;=\; \big(1+\tfrac{\veps}{2}\big)\log t
  \]
  for sufficiently large $t$, as long as $t\ll s\ll t^2/\log t$.

  Hence, for $t$ large enough, $\prob{L_r > 0} \,\leqs\, t\+e^{-(1+\veps/2)\log t} \,=\, t^{-\veps/2}$.
\end{proof}

\begin{figure}[t] 
  \centering
  \begin{tikzpicture}[scale=0.07]
    \fill[gray!20!white] (10.5,.5) -- (10.5,10.5) -- (50.5,50.5) -- (50.5,.5) -- (10.5,.5);
    \foreach \y [count=\x] in {0, 0, 0, 0, 0, 0, 0, 0, 0, 0,
                               1, 4, 0, 0, 6, 0, 5, 3, 0, 3, 3, 2, 2, 6, 3, 2, 9, 6, 1, 3,
                               7, 0, 3, 1, 0, 0, 2, 4, 0, 8, 1, 4, 1, 0, 0, 4, 0, 4, 1, 1}{
      \ifnum0=\y {} \else {
        \foreach \z in {1,...,\y} \fill[radius=0.35,blue!50!red!50!black] (\x,\z) circle;
      } \fi
     \draw (10.5,.5) -- (10.5,10.5) -- (50.5,10.5) -- (50.5,.5);
    }
    \draw[<->] (8.5,.5) -- (8.5,10.5);
    \node[left] at (8.5,5.5) {$r$};
    \draw[<->] (10.5,-1.5) -- (50.5,-1.5);
    \node[below] at (30.5,-1.5) {$t$};
    \node[below] at (30.5,50) {$\Cn[t,s,r]$};
\end{tikzpicture}
$\quad\quad\qquad$
\begin{tikzpicture}[scale=0.07]
    \fill[gray!20!white] (10.5,.5) -- (10.5,10.5) -- (50.5,50.5) -- (50.5,.5) -- (10.5,.5);
    \foreach \y [count=\x] in {0, 0, 0, 0, 0, 0, 0, 0, 0, 0,
                               3, 0, 0, 1,12, 0, 1, 2,18, 3, 2, 1, 0, 4, 2, 1, 0, 0, 1, 3,
                               1, 0,16, 1, 0, 0, 2, 4, 0, 3, 1, 1, 1, 0, 0,13, 0, 2, 0, 1}{
      \ifnum0=\y {} \else {
        \foreach \z in {1,...,\y} \fill[radius=0.35,blue!50!red!50!black] (\x,\z) circle;
      } \fi
     \draw (10.5,.5) -- (10.5,10.5) -- (50.5,50.5) -- (50.5,.5);
    }
    \draw[<->] (8.5,.5) -- (8.5,10.5);
    \node[left] at (8.5,5.5) {$r$};
    \draw[<->] (10.5,-1.5) -- (50.5,-1.5);
    \node[below] at (30.5,-1.5) {$t$};
    \node[below] at (30.5,50) {$\In[t,s,r]$};
\end{tikzpicture}
$\quad\quad\qquad$
\begin{tikzpicture}[scale=0.07]
    \fill[gray!20!white] (10.5,.5) -- (10.5,10.5) -- (50.5,50.5) -- (50.5,.5) -- (10.5,.5);
    \foreach \y [count=\x] in {0, 0, 0, 0, 0, 0, 0, 0, 0, 0,
                               0, 1,21, 0, 0,19, 0, 3, 0, 3, 0, 1, 1, 2, 2, 1,11, 0, 0, 0,
                               1, 0, 2, 1, 0, 0, 2,14, 0, 1, 1, 2, 3, 0, 0, 2, 0, 1, 1, 4}{
      \ifnum0=\y {} \else {
        \foreach \z in {1,...,\y} \fill[radius=0.35,blue!50!red!50!black] (\x,\z) circle;
      } \fi
     \draw (10.5,.5) -- (10.5,50.5);
     \draw (50.5,.5) -- (50.5,50.5);
    }
    \draw[<->] (8.5,.5) -- (8.5,10.5);
    \node[left] at (8.5,5.5) {$r$};
    \draw[<->] (10.5,-1.5) -- (50.5,-1.5);
    \node[below] at (30.5,-1.5) {$t$};
    \node[below] at (30.5,50) {$\Cn[t,s]$};
\end{tikzpicture}

  \caption{A restricted weak composition, an inversion sequence suffix, and an unrestricted weak composition}
  \label{figApproxInvSeqs}
\end{figure}
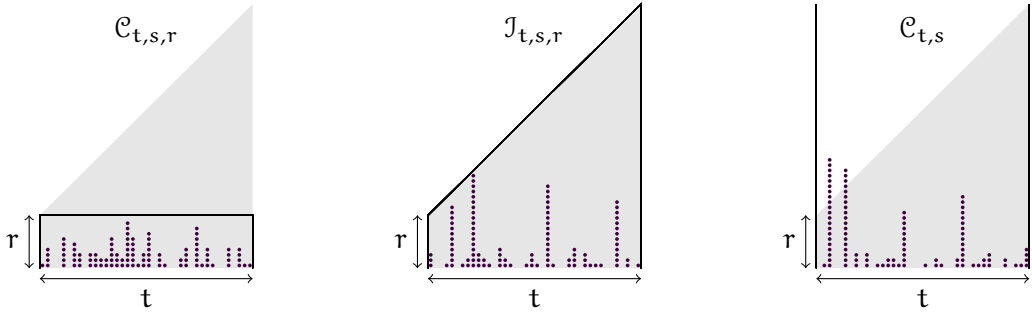

We use weak compositions to approximate all but the first few terms in inversion sequences.
Let $\In[t,s,r]$ denote the set of \emph{inversion sequence suffixes} that consists of weak $t$-compositions $(e_j)_{j=1}^t$ of $s$ in which $e_j<j+r$ for each $j$.
See the middle of Figure~\ref{figApproxInvSeqs} for an illustration.

If $r$ is sufficiently large, then the number of these inversion sequence suffixes may be approximated by the number of unrestricted weak compositions.

\begin{corO}\label{corWeakComp}
  Suppose $t\ll s\ll t^2/\log t$ and $r\geqs(1+\veps)\tfrac{s}{t}\log t$ for some positive $\veps$.
  Then
  \[
  |\In[t,s,r]| \;\sim\; \tbinom{s+t-1}{s}
  .
  \]
\end{corO}
\begin{proof}
Let $\Cn[t,s,r]$ be the set
of \emph{restricted} weak $t$-compositions of $s$,
in which every term is less than~$r$.
Clearly $\Cn[t,s,r]\subset\In[t,s,r]\subset\Cn[t,s]$ (see Figure~\ref{figApproxInvSeqs}).
By Proposition~\ref{propWeakComp}, we have $|\Cn[t,s,r]|\sim|\Cn[t,s]|$ under the specified conditions on $t$, $s$ and $r$.
So, $|\In[t,s,r]|\sim|\Cn[t,s]|=\binom{s+t-1}{s}$.
\end{proof}

\subsection{Counting permutations}

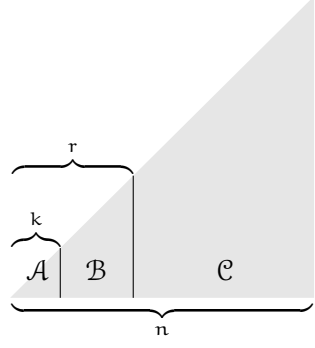
\begin{figure}[t] 
  \centering
\vspace{54pt}
{\small{$\overbrace{\color{white}..................}^r{\color{white}.........................}$}}
\vspace{-69pt}

\vspace{79.5pt}
{\small{$\overbrace{\color{white}.....}^k{\color{white}....................................}$}}
\vspace{-97.5pt}

\begin{tikzpicture}[scale=0.08]
    \fill[gray!20!white] (.5,.5) -- (50.5,50.5) -- (50.5,.5) -- (.5,.5);
     \draw (8.75,.5) -- (8.75,8.75);
     \draw (20.75,.5) -- (20.75,20.75);
\end{tikzpicture}

\vspace{-20.5pt}
$\AAA$ \hspace{8pt} $\BBB$ \hspace{35pt} $\CCC$ \hspace{21pt}

\vspace{-5.75pt}
{\small{$\underbrace{\color{white}.............................................}_n$}}

  \caption{The partitioning of inversion sequences}
  \label{figSplit3}
\end{figure}

To make use of this approximation, we partition the terms of the inversion sequence of an $n$-permutation into three parts.
Given some $k>0$ and $r\geqs k$, part $\AAA$ consists of the first $k$ terms of the sequence (the first $k$ boxes),
$\BBB$ consists of the next $r-k$ terms, and
$\CCC$ consists of the remaining $n-r$ terms.
See Figure~\ref{figSplit3} for an illustration.

We use this tripartition as follows:
Firstly, we place a specific pattern of length $k$ in part $\AAA$.
Secondly, the value of $r$ is chosen so that we can approximate the number of ways of filling part $\CCC$ by using Corollary~\ref{corWeakComp}.
Finally, we sum over each possible way of placing balls in part~$\BBB$.

Suppose we decide to place exactly $\ell$ balls in the first $k$ boxes (part $\AAA$) in some particular way.
Let $\Nknml=|\In[n-k,m-\ell,k]|$
be the number of ways of distributing an additional $m-\ell$ balls among the boxes in parts $\BBB$ and $\CCC$.

Equivalently, $\Nknml$ is the number of $n$-permutations with $m$ inversions whose first $k$ points form some particular permutation having $\ell$ inversions.
For example, $\Nknml[0]$ is the number of $n$-permutations with $m$ inversions whose first $k$ points are increasing.
These have inversion sequences that begin with $k$ zeros (the first $k$ boxes are empty).

If $B=\binom{r}{2}-\binom{k}{2}$ denotes the total capacity of part $\BBB$ (boxes $k+1,\ldots,r$), and
$b_i=|\In[r-k,i,k]|$ is the number of distinct ways of placing exactly $i$ balls in these $r-k$ boxes (for $i=0,\ldots,B$), then
we can express $\Nknml$ as follows:
\begin{equation}\label{eqNkml}
\Nknml \;=\; \sum_{i=0}^B b_i \+ \big|\In[n-r,m-\ell-i,r]\big|
,
\end{equation}
where we sum over the possible choices for the contents of part $\BBB$.

We now approximate the terms in this sum by using Corollary~\ref{corWeakComp}.

\together4
\begin{propO}\label{propNknml}
  If $n\ll m \ll n^2/\log^2n$ and $r=\ceil{\frac{2 m}{n}\log n}$, then
  \[
  \Nknml \;\sim\; 
  \sum_{i=0}^B b_i \+ 
  \binom{m-i \:+\: n-r-1 \:-\: \ell}{n-r-1}
  ,
  \]
  where $B=\binom{r}{2}-\binom{k}{2}$ and $b_i=|\In[r-k,i,k]|$.
\end{propO}

\begin{proof}
To use
Corollary~\ref{corWeakComp}
to approximate $|\In[n-r,m-\ell-i,r]|$
for any nonnegative $\ell\leqs\binom{k}{2}$ and $i\leqs B$, we require the following three inequalities to hold for some
$\veps>0$:
\[
r \;\geqs\; (1+\veps)\frac{m-(\ell+i)}{n-r}\log (n-r)
\text{~~~~~~~and~~~~~~~}
n-r \;\ll\; m-(\ell+i) \;\ll\; \frac{(n-r)^2}{\log (n-r)}
.
\]

Since $m\ll n^2/\log n$, we have $r=\ceil{\frac{2 m}{n}\log n}\ll n$, and thus the first inequality is satisfied for sufficiently large $n$.

Similarly, the third inequality holds as a consequence of $m\ll n^2/\log n$ and $r\ll n$.

Finally, given that $m\ll n^2/\log^2n$, we have
\[
r^2
\;=\; \ceil{\frac{2 m}{n}\log n}^{\!2}
\;\ll\; \frac{2 m \log n}{n} \, \frac{2 n}{\log n}
\;=\; 4 m
.
\]
Thus, $m\gg \binom{r}{2}$, which is the maximum possible value of $\ell+i$.
Together with $n\ll m$, this is sufficient to ensure that the second inequality is satisfied.

\textbf{\color{red} As pointed out by a referee, it is not sufficient, for an unbounded sum, to determine asymptotes of the individual terms. The following statement is thus false.\label{ref1}} 

The result then follows from~\eqref{eqNkml}.
\end{proof}

Note that the proof of this proposition requires $m\ll n^2/\log^2n$.
For faster-growing $m$, one or more of the three inequalities fails to hold, so Corollary~\ref{corWeakComp} cannot be applied.


Our final goal in this section is to compare values of $\Nknml$ for different ranges of values for~$\ell$. We make use of the following two simple results.

\begin{propO}\label{propBinomialRatio}
  If $ x\ll  y$,
  then
  \[
  \liminfty[ x]
  \binom{ y}{ x}{\Big/}\binom{ y-\delta}{ x} \;=\;
  \begin{cases}
  1           & \text{if~ $\delta\ll  y/ x$,} \\
  e^{\alpha}  & \text{if~ $\delta\sim\alpha\+  y/ x$, ~for any~ $\alpha>0$,} \\
  \infty      & \text{if~ $\delta\gg  y/ x$.}
  \end{cases}
  \]
\end{propO}

\begin{proof}
  By Stirling's approximation,
  \[
  \tbinom{ y}{ x}{\Big/}\tbinom{ y-\delta}{ x}
      \;\sim\; \Big(1-\tfrac{ x\delta}{ (y-x)( y-\delta)}\Big)^{\! y+\half} \,
               \Big(1+\tfrac{\delta}{ y- x-\delta}\Big)^{\! x} \,
               \Big(1+\tfrac{ x}{ y-x-\delta}\Big)^{\!\delta}
  ,
  \]
  from which the result can be seen to follow.
\end{proof}

\begin{propO}\label{propSumBounds}
  Suppose that we have positive $a_i,x_i,y_i$ for $i=1,\ldots,n$ and that there are $L,U$ such that $L\leqs x_i/y_i\leqs U$ for each $i$.
  If $X=\sum_{i=1}^na_ix_i$ and $Y=\sum_{i=1}^na_iy_i$,
  then $L\leqs X/Y\leqs U$.
\end{propO}

\begin{proof}
  \[
  Y\+L \;=\;
  \sum_{i=1}^na_iy_i L \;\leqs\;
  \sum_{i=1}^na_ix_i \;=\;
  X \;=\;
  \sum_{i=1}^na_ix_i \;\leqs\;
  \sum_{i=1}^na_iy_i U \;=\;
  Y\+U
  .
  \qedhere
  \]
\end{proof}

Recall that $\Nknml$ is the number of $n$-permutations with $m$ inversions whose first $k$ points form some particular permutation having $\ell$ inversions.
The following proposition establishes the threshold for change in the asymptotic value of $\Nknml$ for semi-sparse permutations in terms of the growth of $\ell$ with respect to $n$ and $m$.

\begin{propO}\label{propNComp}
  Suppose $n\ll m \ll n^2/\log^2n$. Then,
  \begin{align*}
    \Nknml &\;\sim\; \Nknml[0]               && \hspace{-3cm} \text{if~ $\ell \ll m/n$,} \\[2pt]
    \Nknml &\;\sim\; e^{-\beta}\Nknml[0]     && \hspace{-3cm} \text{if~ $\ell \sim \beta m/n$,} \\[2pt]
    \Nknml &\;\ll\; \Nknml[0]                && \hspace{-3cm} \text{if~ $\ell \gg m/n$.}
  \end{align*}
\end{propO}

\begin{proof}
Suppose $r\ll n$ and $i\ll m$.
Then $\frac{m-i + n-r-1}{n-r-1} \sim \frac{m}{n}$.
Set $r=\ceil{\frac{2 m}{n}\log n}$.

Suppose $\ell\ll m/n$.
Then, by Proposition~\ref{propBinomialRatio}, we have $\binom{m-i + n-r-1 \:-\: \ell}{n-r-1} \sim \binom{m-i + n-r-1}{n-r-1}$, and so, by Propositions~\ref{propNknml} and~\ref{propSumBounds}, $\Nknml\sim \Nknml[0]$.

Similarly, if $\ell\sim \beta m / n$ then $\binom{m-i + n-r-1 \:-\: \ell}{n-r-1} \sim e^{-\beta}\binom{m-i + n-r-1}{n-r-1}$, and thus $\Nknml\sim e^{-\beta}\Nknml[0]$.

Finally, suppose $\ell\gg m/n$.
Then, by Proposition~\ref{propBinomialRatio},
$\binom{m-i + n-r-1 \:-\: \ell}{n-r-1} \;\ll\; \binom{m-i + n-r-1}{n-r-1}$.
Hence, by Propositions~\ref{propNknml} and~\ref{propSumBounds}, $\Nknml[\ell]\ll \Nknml[0]$.
\end{proof}

\section{Threshold for local uniformity}\label{sectUnif}

We are now in a position to establish the threshold for local uniformity. 
First we prove that a semi-sparse permutation $\ssnm$ is indeed locally uniform. This is Theorem~\ref{thmPatUnif}, which we restate here.

\textbf{Theorem~\ref{thmPatUnif}.}
  \emph{
  Suppose $n\ll m \ll n^2/\log^2n$ and $k \ll \sqrt{m/n}$.
  Then, for any sequence of permutations $(\tau_n)$ with
  $|\tau_n|=k$
  and
  any sequence of positive integers $(j_n)$ with $j_n\leqs n+1-k$,
  \[
  \prob{\text{$\tau_n$ occurs at position $j_n$ in $\ssnm$}}
  \;\sim\;
  \frac{1}{k!}
  .
  \]}\vspace{-12pt}

\begin{proof}
  Since $|\tau_n|=k\ll \sqrt{m/n}$,
  we have $\inv(\tau_n)\leqs\binom{k}{2}\ll m/n$.
  Therefore, by Proposition~\ref{propNComp}, we have $\Nknml[\inv(\tau_n)]\sim \Nknml[0]$, and hence
  \[
  \prob{\text{$\ssnm[m][1,k]$ forms $\tau_n$}}
  \;\sim\;
  \prob{\text{$\ssnm[m][1,k]$ is increasing}}
  ,
  \]
  and this is true whatever sequence $(\tau_n)$ is selected.
  Thus, asymptotically, every possible choice for the pattern formed by the first $k$ points of $\ssnm$ is equally probable, and so
  \[
  \prob{\text{$\tau_n$ occurs at position $1$ in $\ssnm$}}
  \;\sim\;
  \frac{1}{k!}
  .
  \]
  The result then follows from the position independence
  of the distribution of consecutive patterns
  in $\ssnm$ (Proposition~\ref{propShiftInvariance}).
\end{proof}

We now consider the behaviour in the critical window.
In order to do this we require tail bounds on the distribution of the number of inversions in a random $n$-permutation.
We use $\ssn$ to denote a permutation chosen uniformly 
from $\Sn$.

\begin{propO}\label{propTailBounds}
  For any $\theta>0$,
  \[
  \prob{\big|\dinv(\ssn) - \thalf\big| \:>\: \theta } \;<\; 2\+e^{-\theta^2n}
  .
  \]
\end{propO}

\begin{proof}
  We apply Hoeffding's inequality, which states that if $X_1,\ldots,X_n$ are independent random variables such that $0\leqs X_i\leqs u_i$ for all $i$, and $S_n=\sum_{i=1}^nX_i$, then
  \[
  \prob{\big|S_n-\bbE[S_n]\big| \:>\: t} \;<\; 2 \exp\!\left( - \+ \frac{2\+t^2}{\sum_{i=1}^n u_i^2} \right).
  \]
  Now, as illustrated by the balls-in-boxes model,
  $\inv(\ssn) \sim \sum_{i=1}^n \mathsf{Unif}[0,i-1]$, the sum of $n$ independent discrete uniform random variables, so we can set $u_i=i-1$, yielding
  \[
  \sum_{i=1}^n u_i^2
  \;=\; \frac{(2n-1)n(n-1)}{6}
  ,
  \]
  which does not exceed 
  $\frac{2}{n}\binom{n}{2}^2$, from which the result follows directly.
\end{proof}

In the critical window, where $k\sim\alpha\sqrt{m/n}$, the asymptotic probability of a particular consecutive pattern depends (only) on its inversion density.

\textbf{Theorem~\ref{thmPatCrit}.} \emph{
  Suppose $n\ll m \ll n^2/\log^2n$ and $k\sim \alpha\sqrt{m/n}$ for some $\alpha>0$.
  Fix $\rho\in[0,1]$.
  Then, for any sequence of permutations $(\tau_n)$ with $|\tau_n|=k$ and
  $\dinv(\tau_n)\sim\rho$,
  and
  any sequence of positive integers $(j_n)$ with $j_n\leqs n+1-k$,
  \[
  \prob{\text{$\tau_n$ occurs at position $j_n$ in $\ssnm$}}
  \;\sim\;
  { e^{(1-2\rho)\alpha^2/4} }
  \frac {1}
  { k! }
  .
  \]}\vspace{-12pt}

\begin{proof}

Since $|\tau_n|=k\sim\alpha\sqrt{m/n}$ and $\inv(\tau_n)\sim\rho \binom{k}{2}$,
we have $\inv(\tau_n)\sim\frac{\rho\alpha^2m}{2n}$.
Therefore, by Proposition~\ref{propNComp}, we have $\Nknml[\inv(\tau_n)]\sim e^{-\rho\alpha^2/2}\Nknml[0]$, and hence
\begin{equation}\label{eqAsymRatio}
  \prob{\text{$\ssnm[m][1,k]$ forms $\tau_n$}}
  \;\sim\;
  e^{-\rho\alpha^2/2} \, \prob{\text{$\ssnm[m][1,k]$ is increasing}}
  .
\end{equation}

Now clearly, for any given $k$, we have 
  $
  \sum_{\vphi\in[0,1]} \prob{ \dinv(\ssnm[m][1,k]) = \vphi}
  =
  1 ,
  $
  where the sum should be understood to be
  over the finite set of possible inversion densities of $k$-permutations.
  Equivalently,
  \[
  \sum_{\vphi\in[0,1]} \Big|\Sn[k,\vphi\binom{k}{2}]\Big| \,  \prob{\text{$\ssnm[m][1,k]$ forms $\pi_k^\vphi$}}  \;=\;  1 ,
  \]
  where, for each valid value of $\vphi$, we choose $\pi_k^\vphi$ to be some $k$-permutation with inversion density exactly $\vphi$.

  There are exactly $k!\,\prob{\dinv(\ssk)=\vphi}$ permutations of length $k$ with inversion density $\vphi$.
  So,
  \[
  k! \sum_{\vphi\in[0,1]} \prob{\dinv(\ssk)=\vphi} \,  \prob{\text{$\ssnm[m][1,k]$ forms $\pi_k^\vphi$}}  \;=\;  1 .
  \]
  Now, in an analogous manner to~\eqref{eqAsymRatio}, for every valid $\vphi$ we have
  \[
  \prob{\text{$\ssnm[m][1,k]$ forms $\pi_k^\vphi$}}
  \;\sim\;
  e^{-\vphi\alpha^2/2} \, \prob{\text{$\ssnm[m][1,k]$ is increasing}}
  ,
  \]
  where we take limits over those $n$ for which $\vphi \binom{k}{2}\in\bbN$.
  Thus,
\begin{equation}\label{eqAsymSum1}
  k! \, \prob{\text{$\ssnm[m][1,k]$ is increasing}}
  \sum_{\vphi\in[0,1]}
  \prob{\dinv(\ssk)=\vphi}
  e^{-\vphi\alpha^2/2}
  \;\sim\;
  1
  .
\end{equation}

  We now make use of our tail bounds for the inversion density. By Proposition~\ref{propTailBounds},
  \[
  \prob{ \big| \dinv(\ssk) - \thalf \big| \:>\: k^{-1/4} }
  \;<\;
  2\+e^{-\sqrt{k}}
  ,
  \]
  and $e^{-\vphi\alpha^2/2}$ is no greater than $1$ for any $\vphi$.
  So the contribution to the sum in~\eqref{eqAsymSum1} from inversion densities
  that differ from $\nhalf$ by
  more than $k^{-1/4}$ is less than $2\+e^{-\sqrt{k}}$, which tends to zero.

  On the other hand, the remaining contribution to the sum (from inversion densities close to $\nhalf$) lies between
  \[
  \big(1-2\+e^{-\sqrt{k}}\big) e^{-(1+k^{-1/4})\alpha^2/4}
  \text{~~~~~and~~~~~}
  e^{-(1-k^{-1/4})\alpha^2/4}
  ,
  \]
  and so tends to $e^{-\alpha^2/4}$.
  Thus,
  $
  k! \, \prob{\text{$\ssnm[m][1,k]$ is increasing}} 
  \sim
  e^{\alpha^2/4}
  .
  $

  Applying~\eqref{eqAsymRatio} then yields
  \[
  \prob{\text{$\tau_n$ occurs at position $1$ in $\ssnm$}}
  \;\sim\;
  { e^{(1-2\rho)\alpha^2/4} }
  \frac {1}
  { k! }
  ,
  \]
  and the result follows from the position independence
  of the distribution of consecutive patterns
  in $\ssnm$ (Proposition~\ref{propShiftInvariance}).
\end{proof}

Any larger subpermutation with a sufficient number of inversions almost never occurs.

\textbf{Theorem~\ref{thmPat0}.} \emph{
  Suppose $n\ll m \ll n^2/\log^2n$ and $k\gg \sqrt{m/n}$.
  Suppose $m/nk^2\ll\rho\leqs1$.
  Then, for any sequence of permutations $(\tau_n)$ with $|\tau_n|=k$ and
  $\dinv(\tau_n)\sim\rho$,
  and
  any sequence of positive integers $(j_n)$ with $j_n\leqs n+1-k$,
  \[
  \prob{\text{$\tau_n$ occurs at position $j_n$ in $\ssnm$}}
  \;\ll\;
  \prob{\text{$\ssnm[m][j_n,j_n+k-1]$ is increasing}}
  .
  \]}\vspace{-12pt}

\begin{proof}
Since $|\tau_n|=k\gg\sqrt{m/n}$ and $\inv(\tau_n)\sim\rho \binom{k}{2}$,
with $\rho\gg m/nk^2$,
we have $\inv(\tau_n)\gg m/n$.
Therefore, by Proposition~\ref{propNComp}, we have $\Nknml[\inv(\tau_n)]\ll \Nknml[0]$, and hence
  \[
  \prob{\text{$\ssnm[m][1,k]$ forms $\tau_n$}}
  \;\ll\;
  \prob{\text{$\ssnm[m][1,k]$ is increasing}}
  .
  \]
The result then follows from the position independence
  of the distribution of consecutive patterns
  in $\ssnm$ (Proposition~\ref{propShiftInvariance}).
\end{proof}

Theorems~\ref{thmPatUnif}, \ref{thmPatCrit} and~\ref{thmPat0} thus reveal
the existence of a threshold at $k=\sqrt{m/n}$ for
the uniform distribution of
consecutive $k$-subpermutations of semi-sparse $\ssnm$.

\begin{figure}[ht]
  \centering
\begin{tikzpicture}[scale=0.335]
    \foreach \y [count=\x] in {0,1,2,2,0,4,2}{
      \node[blue!50!black] at (\x,-1) {\small\textsf{\y}};
      \draw[very thin] (\x+1,.48)--(\x,-.5);
    }
    \node[red!50!black] at (8,-1) {\small\textsf{2}};
    \node[green!50!black] at (9,-1) {\small\textsf{3}};
    \draw[very thin] (1,.5)--(1,.25)--(8,-.25)--(8,-.5);
    \draw[very thin] (9,.5)--(9,-.5);
    \plotpermgrid  [blue!50!black]{9}{0,8,5,1,4,9,2,7,0}
    \plotpermnobox[red!50!black]{9}{3}
    \plotpermnobox[green!50!black]{9}{0,0,0,0,0,0,0,0,6}
    \circpt{1}{3}
    \circpt{9}{6}
\end{tikzpicture}
$\qquad\qquad\qquad\quad$
\begin{tikzpicture}[scale=0.335]
    \foreach \y [count=\x] in {0,1,2,2,0,4,2}{
      \node[blue!50!black] at (\x,-1) {\small\textsf{\y}};
      \draw[very thin] (\x+1,.48)--(\x,-.5);
    }
    \node[red!50!black] at (8,-1) {\small\textsf{2}};
    \node[green!50!black] at (9,-1) {\small\textsf{7}};
    \draw[very thin] (1,.5)--(1,.25)--(8,-.25)--(8,-.5);
    \draw[very thin] (9,.5)--(9,-.5);
    \plotpermgrid  [blue!50!black]{9}{0,8,6,1,5,9,3,7,0}
    \plotpermnobox[red!50!black]{9}{4}
    \plotpermnobox[green!50!black]{9}{0,0,0,0,0,0,0,0,2}
    \circpt{1}{4}
    \circpt{9}{2}
\end{tikzpicture}
  \caption{Permutations built by adjoining new initial and final points, together with their nonstandard inversion sequences}
  \label{figAddInversions}
\end{figure}
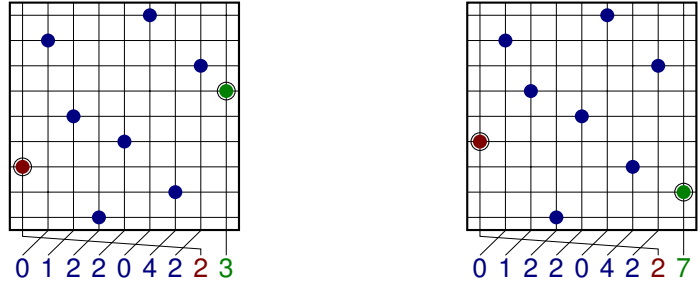

\section{Threshold for inversions}\label{sectInv}

We now turn our attention to the uniformity of \emph{inversions}: How close do indices $i$ and $j$ need to be for
$i,j$ to be as likely to form an inversion in $\ssnm$ as not?
We begin with some counting.

Suppose $\pi$ is a $(k-1)$-permutation.
Let us consider how new first and last points may be adjoined to $\pi$ so that the resulting permutation has exactly $\ell$ more inversions than $\pi$.
Specifically, we want to determine how many distinct ways there are to construct a $(k+1)$-permutation $\tau$ such that $\tau[2,k]$ forms $\pi$ and $\inv(\tau)-\inv(\pi)=\ell$.
The answer depends on $k$ and $\ell$, and on whether or not $\tau(1)\tau(k+1)$ forms an inversion.
It doesn't depend on any other properties of $\pi$ or~$\tau$.

Suppose $0 \leqs \ell\leqs k-1$.
In this case, $\inv(\tau)-\inv(\pi)=\ell$ precisely when $\tau(k+1)-\tau(1)=k-\ell$.
For example, on the left of Figure~\ref{figAddInversions}, we have $\tau(9)-\tau(1)=3$ with $k=8$ and $\ell=5$.
Note that $\tau(1)\tau(k+1)$ does not form an inversion.
In this case it is readily checked that there are exactly $\ell+1$ distinct ways to construct an appropriate $\tau$ from any given~$\pi$, adding $\ell$ inversions.

On the other hand, if $k\leqs \ell\leqs2\+k-1$, then $\inv(\tau)-\inv(\pi)=\ell$ whenever $\tau(1)-\tau(k+1)=\ell+1-k$.
For example, on the right of Figure~\ref{figAddInversions}, we have $\tau(1)-\tau(9)=2$ with $k=8$ and $\ell=9$.
Note that, in this case, $\tau(1)\tau(k+1)$ does form an inversion. 
It can be seen that there are exactly $2k-\ell$ distinct ways to build a suitable $\tau$ from any given $(k-1)$-permutation~$\pi$, adding $\ell$ inversions.

\together3
We use nonstandard inversion sequences to represent permutations built this way, reflecting how they are constructed. 
Specifically, we represent a $(k+1)$-permutation $\tau$ constructed from a $(k-1)$-permutation $\pi$
by a sequence $(e_j)_{j=1}^{k+1}$, where the first ${k-1}$ terms 
form the standard inversion sequence for $\pi$, so
\[
e_j \;=\; \big|\{i: 1\leqs i<j \text{~and~} \tau(i+1)>\tau(j+1) \}\big|
\]
for $1\leqs j\leqs k-1$, and the final two terms are given by
\begin{align*}
e_k     & \;=\; \big|\{i: 2\leqs i\leq k \text{~and~} \tau(1)>\tau(i) \}\big| , \\[2pt]
e_{k+1} & \;=\; \big|\{i: 1\leqs i<k+1 \text{~and~} \tau(i)>\tau(k+1) \}\big| ,
\end{align*}
recording the number of inversions created by adjoining the new initial and final points, respectively.
See Figure~\ref{figAddInversions} for two examples.
Note that we still have $e_j<j$ for each $j$.

\begin{figure}[ht]
  \centering
\vspace{54pt}
{\small{$\overbrace{\color{white}..................}^r{\color{white}.........................}$}}
\vspace{-69pt}

\vspace{79.5pt}
{\small{$\overbrace{\color{white}.....}^k{\color{white}....................................}$}}
\vspace{-97.5pt}

\begin{tikzpicture}[scale=0.08]
    \fill[gray!20!white] (.5,.5) -- (50.5,50.5) -- (50.5,.5) -- (.5,.5);
     \draw (8.75,.5) -- (8.75,8.75);
     \draw (20.75,.5) -- (20.75,20.75);
\end{tikzpicture}

\vspace{-20.5pt}
$\AAA$ \hspace{8pt} $\BBB$ \hspace{35pt} $\CCC$ \hspace{21pt}

\vspace{-5.75pt}
{\small{$\underbrace{\color{white}.............................................}_n$}}

  \caption{The partitioning of modified inversion sequences}
  \label{figSplit3mod}
\end{figure}
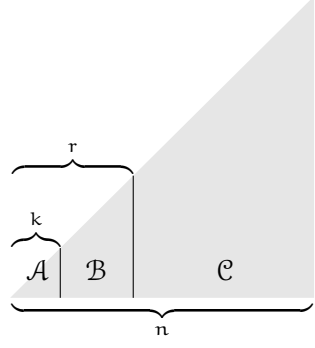

We now construct a modified inversion sequence for each $n$-permutation $\sigma$, in which the first $k+1$ terms form the nonstandard inversion sequence for $\sigma[1,k+1]$ and subsequent terms are standard.
Thus, the $k$th and $(k+1)$th terms record the number of inversions created by adjoining $\sigma(1)$ and $\sigma(k+1)$ to $\sigma[2,k]$.

As before, for a given $r\geqs k$, we partition the terms (or boxes) of modified inversion sequences into three parts (see Figure~\ref{figSplit3mod}).
Let $A=\binom{r}{2}-(2k-1)$ be the total capacity of
parts $\AAA$ and $\BBB$ excluding the ``special'' boxes $k$ and $k+1$.
That is, $A$ is the capacity of boxes $1,\ldots,k-1$ and $k+2,\ldots,r$.
Now, for each $i=0,\ldots,A$,
let
$a_i$ be the number of distinct ways of placing exactly $i$ balls in these $r-2$ boxes.

Then the number of $n$-permutations $\sigma$ with $m$ inversions in which $\sigma(1)<\sigma(k+1)$ is given by
\begin{equation}\label{eqNknmup}
\Nknmup
\;=\;
\sum_{i=0}^A a_i
\sum_{\ell=0}^{k-1} \+
(\ell+1) \+
\big|\In[n-r,m-\ell-i,r]\big|
,
\end{equation}
where $\ell+1$ is the number of ways of choosing $\sigma(1)$ and $\sigma(k+1)$, or equivalently the contents of boxes $k$ and $k+1$, so as to contribute $\ell$ inversions.

Similarly, the number of $n$-permutations $\sigma$ with $m$ inversions in which $\sigma(1)$ and $\sigma(k+1)$ form an inversion is given by
\begin{align}
\Nknmdn
& \;=\;
\sum_{i=0}^A a_i
\sum_{\ell=k}^{2k-1} \+
(2k-\ell) \+
\big|\In[n-r,m-\ell-i,r]\big| \nonumber\\[2pt]
& \;=\;
\sum_{i=0}^A a_i
\sum_{\ell=0}^{k-1} \+
(\ell+1) \+
\big|\In[n-r,m-\ell-i-(2k-2\ell-1),r]\big| \label{eqNknmdn}
,
\end{align}
where the second expression results from the change of variable $\ell=2k-1-\ell$.

As we did in Proposition~\ref{propNknml}, we now use Corollary~\ref{corWeakComp} to approximate the terms of these sums.

\begin{propO}\label{propNdnmx}
  If $n\ll m \ll n^2/\log^2n$ and $r=\ceil{\frac{2 m}{n}\log n}$, then
  \begin{align*}
  \Nknmup
  & \;\sim\;
  \sum_{i=0}^A a_i
  \sum_{\ell=0}^{k-1} \+
  (\ell+1) \+
  \binom{m-\ell-i \:+\: n-r-1}{n-r-1}
  , \\[2pt]
  \Nknmdn
  & \;\sim\;
  \sum_{i=0}^A a_i
  \sum_{\ell=0}^{k-1} \+
  (\ell+1) \+
  \binom{m-\ell-i \:+\: n-r-1 \:-\: (2k-2\ell-1)}{n-r-1}
  ,
  \end{align*}
  where $A=\binom{r}{2}-(2k-1)$ and $a_i$ is as defined above.
\end{propO}

\begin{proof}

Reasoning almost identical to that in Proposition~\ref{propNknml} shows that the choice of $r$ and the bounds on $n$ and $m$ are sufficient
to ensure that Corollary~\ref{corWeakComp} can be used
to approximate $|\In[n-r,m-\ell-i,r]|$ and $|\In[n-r,m-\ell-i-(2k-2\ell-1),r]|$
for any nonnegative
$i\leqs A$ and $\ell<k$.
The result then follows from~\eqref{eqNknmup} and~\eqref{eqNknmdn}.
\end{proof}

Using these approximations, we can now exhibit a threshold at $k=m/n$ for the uniformity of {inversions}.

\textbf{Theorem~\ref{thmInvUnif}.} \emph{
  Suppose $n\ll m \ll n^2/\log^2n$ and $(j_n)$ is any sequence of positive integers such that $j_n\leqs n-k$.
  Then,
  \[
  \liminfty\prob{\ssnm(j_n) > \ssnm(j_n+k)}
  \;=\;
  \begin{cases}
  \thalf & \text{if $k\ll {m/n}$,} \\[2pt]
  0 & \text{if $k\gg {m/n}$.}
  \end{cases}
  \]}\vspace{-12pt}

\begin{proof}
  We compare $\Nknmdn$ against $\Nknmup$.

  As long as $r\ll n$ and $\ell<k$,
  by Proposition~\ref{propBinomialRatio},
  \[
  \liminfty
  \tbinom{m-\ell-i \:+\: n-r-1 \:-\: (2k-2\ell-1)}{n-r-1}
  \big/
  \tbinom{m-\ell-i \:+\: n-r-1}{n-r-1}
  \;=\;
  \begin{cases}
  1 & \text{if $k\ll {m/n}$,} \\[2pt]
  0 & \text{if $k\gg {m/n}$.}
  \end{cases}
  \]
  \textbf{\color{red}As pointed out by a referee, this only holds if $k\gg\ell$.\label{ref2}}

  Set $r=\ceil{\frac{2 m}{n}\log n}$.
  If $k\ll {m/n}$, then by Propositions~\ref{propNdnmx} and~\ref{propSumBounds}, we have
  $\Nknmdn\sim \Nknmup$, and so
  \[
  \prob{\text{$\ssnm(1)>\ssnm(k+1)$}}
  \;\sim\;
  \prob{\text{$\ssnm(1)<\ssnm(k+1)$}}
  \;\sim\;
  \thalf
  .
  \]

  On the other hand, if $k\gg {m/n}$, then $\Nknmdn\ll \Nknmup$, and so
  $\prob{\text{$\ssnm(1)<\ssnm(k+1)$}} \sim 1.$

  The result then follows from the position independence
  of the distribution of consecutive patterns
  in $\ssnm$ (Proposition~\ref{propShiftInvariance}).
\end{proof}

In the {critical window}, we have the following behaviour.

\textbf{Theorem~\ref{thmInvCrit}.} \emph{
  Suppose $n\ll m \ll n^2/\log^2n$ and $k\sim \alpha {m/n}$ for some $\alpha>0$.
  Then, for any sequence of positive integers $(j_n)$ with $j_n\leqs n-k$,
  \[
  \liminfty\prob{\ssnm(j_n) > \ssnm(j_n+k)}
  \;=\;
  \frac{e^\alpha(\alpha-1)+1}{(e^\alpha-1)^2}
  \;<\;
  \half
  .
  \]}\vspace{-12pt}

\begin{proof}
It can be checked (using a computer algebra package or otherwise) that
\[
S_{{}^\nearrow}
\;=\;
\sum_{\ell=0}^{k-1} \+ (\ell+1) \+ \binom{y-\ell}{x}
\;=\;
\frac{
(y+1)(y+2)
\binom{y}{x} - (y+1-k) (y+2+x\+k +k) \binom{y-k}{x}}{(x+1) (x+2)}
,
\]
and
\begin{align*}
S_{{}^\searrow}
& \;=\;
\sum_{\ell=0}^{k-1} \+
  (\ell+1) \+
  \binom{y-\ell-(2k - 2\ell - 1)}{x} \\[2pt]
& \;=\;
\frac{(y-x-2 k)(y+1-x-2 k)  \binom{y+1-2 k}{x} - (y+1-x-k) (y-x k -x-3 k) \binom{y+1-k}{x}}{(x+1) (x+2)}
.
\end{align*}

Thus, if $k\sim \alpha\+ y/x$ and $x\ll y$, then by Proposition~\ref{propBinomialRatio},
\begin{align*}
S_{{}^\nearrow}
& \;\sim\;
\frac{y^2}{x^2}\binom{y}{x}
\big(1-(1+\alpha)e^{-\alpha}\big) , \\[2pt]
S_{{}^\searrow}
& \;\sim\;
\frac{y^2}{x^2}\binom{y}{x}
\big(e^{-2\alpha}-(1-\alpha)e^{-\alpha}\big).
\end{align*}

So, if we let $x=n-r-1$ and $y=m-i+x$, by Propositions~\ref{propNdnmx} and~\ref{propSumBounds}
we have
  the following for $p_\alpha=\prob{\ssnm(1) > \ssnm(k+1)}$, the probability that
  $1,k+1$
  forms an inversion in $\ssnm$:
  \[
  p_\alpha
  \;=\;
  \frac{\Nknmdn}{\Nknmup\,+\,\Nknmdn}
  \;\sim\;
  \frac{e^{-2\alpha}-(1-\alpha)e^{-\alpha}}{1-(1+\alpha)e^{-\alpha} \:+\: e^{-2\alpha}-(1-\alpha)e^{-\alpha}
  }
  \;=\;
  \frac{e^\alpha(\alpha-1)+1}{(e^\alpha-1)^2} .
  \]
  For $\alpha>0$, this probability decreases as $\alpha$ increases, since
  \[
  \frac{d\+p_\alpha}{d\alpha}
  \;=\;
  -\frac{e^{\alpha } \left(e^{\alpha } (\alpha -2)+\alpha +2\right)}{(e^{\alpha }-1)^3}
  \]
  and
  \[
  e^{\alpha } (\alpha -2)+\alpha +2
  \;=\;
  \sum_{n\geqs3} (n - 2) \alpha^n/n!
  ,
  \]
  each term in the Taylor expansion being positive. Thus, for
  positive $\alpha$, we have $p_\alpha<\half$, since 
  $\lim_{\alpha\to0}p_\alpha=\half$.

  The result then follows from the position independence
  of the distribution of consecutive patterns
  in $\ssnm$ (Proposition~\ref{propShiftInvariance}).
\end{proof}

\subsection*{Acknowledgements}
The
catalyst for
this work was a brief conversation with
Jakub Slia\v{c}an at
the 26th British Combinatorial Conference in July 2017 at the University of Strathclyde in Glasgow.
I am also grateful to Thomas Selig for acting as a sounding-board for some of my initial ideas
on this and related topics.
\emph{Mathematica}~\cite{Mathematica} was used for experimentation and algebraic manipulation.

\emph{Soli Deo gloria!}


\bibliographystyle{plain}
{\footnotesize\bibliography{../bib/mybib}}

\begin{thebibliography}{10}

\bibitem{AP2013}
H{\"u}seyin Acan and Boris Pittel.
\newblock On the connected components of a random permutation graph with a
  given number of edges.
\newblock {\em J.~Combin. Theory Ser.~A}, 120(8):1947--1975, 2013.

\bibitem{BevanPPBasics}
David Bevan.
\newblock Permutation patterns: basic definitions and notation.
\newblock {\em \href{http://arxiv.org/pdf/1506.06673}{arXiv:1506.06673}}, 2015.

\bibitem{Bevan2022}
David Bevan.
\newblock Independence of permutation limits at infinitely many scales.
\newblock {\em J.~Combin. Theory Ser.~A}, 186:~Paper 105557, 19~pp., 2022.

\bibitem{BonaBook}
Mikl{\'o}s B{\'o}na.
\newblock {\em Combinatorics of Permutations}.
\newblock CRC Press, second edition, 2012.

\bibitem{BP2020}
Jacopo Borga and Raul Penaguiao.
\newblock The feasible region for consecutive patterns of permutations is a
  cycle polytope.
\newblock {\em Algebraic Combinatorics}, 3(6):1259--1281, 2020.

\bibitem{DG1977}
Persi Diaconis and R.~L. Graham.
\newblock Spearman's footrule as a measure of disarray.
\newblock {\em J.~Royal Statist. Soc. Ser.~B}, 39:262--268, 1977.

\bibitem{ER1960}
P.~Erd{\H{o}}s and A.~R{\'e}nyi.
\newblock On the evolution of random graphs.
\newblock {\em Magyar Tud. Akad. Mat. Kutat\'o Int. K\H{o}zl.}, 5:17--61, 1960.

\bibitem{KKRW2020}
Richard Kenyon, Daniel Kr\'al', Charles Radin, and Peter Winkler.
\newblock Permutations with fixed pattern densities.
\newblock {\em Random Structures \& Algorithms}, 56(1):220--250, 2020.

\bibitem{KitaevBook}
Sergey Kitaev.
\newblock {\em Patterns in Permutations and Words}.
\newblock Springer, 2011.

\bibitem{Knuth1973}
Donald~E. Knuth.
\newblock {\em The Art of Computer Programming. {V}ol. 3: {S}orting and
  Searching}.
\newblock Addison-Wesley, 1973.

\bibitem{VatterSurvey}
Vincent Vatter.
\newblock Permutation classes.
\newblock In Mikl{\'o}s B{\'o}na, editor, {\em The Handbook of Enumerative
  Combinatorics}. CRC Press, 2015.

\bibitem{Mathematica}
{Wolfram Research, Inc.}
\newblock Mathematica.
\newblock \href{http://www.wolfram.com/mathematica}{wolfram.com/mathematica}.

\end{thebibliography}

\end{document}